\documentclass[11pt,letterpaper]{article}

\usepackage[T1]{fontenc}
\usepackage{lmodern}
\usepackage[margin=1in]{geometry}
\usepackage{amsmath,amssymb,amsthm}
\usepackage{microtype}
\usepackage[skip=0.75\baselineskip,indent=0pt]{parskip}
\usepackage{xcolor}
\definecolor{linkblue}{RGB}{0,56,164}
\usepackage[colorlinks=true,linkcolor=linkblue,citecolor=linkblue,urlcolor=linkblue,bookmarksnumbered=true]{hyperref}
\hypersetup{pdftitle={Avoiding Isosceles Triangles in the Integer Grid},pdfauthor={Samuel Korsky}}
\newtheorem{theorem}{Theorem}[section]
\newtheorem{lemma}[theorem]{Lemma}
\newcommand{\F}{\mathbb F}
\newcommand{\arxivref}[1]{\href{https://arxiv.org/abs/#1}{arXiv:#1}}
\newcommand{\greenpaper}{\href{https://people.maths.ox.ac.uk/greenbj/papers/open-problems.pdf}{\emph{100 open problems}}}
\newcommand{\numberthis}{\stepcounter{equation}\tag{\theequation}}
\title{Avoiding Isosceles Triangles in the Integer Grid}
\author{Samuel Korsky}
\date{November 18, 2026}

\begin{document}
{\setlength{\parskip}{0pt}\maketitle}

\smallskip

\begin{abstract}
\noindent We construct a set of $\Omega(n)$ points in the $n\times n$ integer grid containing no isosceles triangle (including degenerate triangles). For every $\varepsilon>0$, we also construct a set of $\Omega_\varepsilon(n^{2\log_5 3-\varepsilon})$ grid points containing no isosceles right triangle in any orientation, improving upon a bound of K\'arolyi and Solymosi.
\end{abstract}

\section{Introduction}
\label{sec:intro}

For a positive integer $n$, write $I_n=\{0,1,\ldots,n-1\}$. A set $A\subset I_n^2$ is \emph{isosceles-free} if $\lVert b-a\rVert_2\ne\lVert c-a\rVert_2$ for every three distinct points $a,b,c\in A$. Thus three equally spaced points on a line are forbidden. Let $F(n)$ be the largest cardinality of such a set.

\begin{theorem}\label{thm:main}
There is an absolute constant $c>0$ such that $F(n)\ge cn$ for all $n\ge1$.
\end{theorem}

In 1973, Erd\H{o}s~\cite{erdos} attributed to Riddell the question of how large an isosceles-free subset one can always extract from $m$ arbitrary points in the plane. He singled out the square lattice as a possible extremal example. On page 120 of that paper, he suggested that a method of his and Guy's could yield $cn$ points in the grid, and speculated that for some $\varepsilon>0$ one might even find $n^{1+\varepsilon}$ such points. He also suggested an upper bound with exponent strictly below two. Another paper of Erd\H{o}s and Guy~\cite{erdosguy} concerns the related requirement that \emph{all} pairwise distances be different. In an isosceles-free set a distance may occur repeatedly on disjoint pairs of points.

Wu~\cite{wu} studied the rectangular-grid problem with positive-area triangles. Ellenberg~\cite{ellenberg} explicitly asked the square-grid question including collinear triples. Motivated by the problem of reconstructing point configurations from comparisons of distances, Ellenberg and Jain~\cite{ellenbergjain} asked whether, for small $\varepsilon>0$, the $n\times n$ grid contains an isosceles-free set of size $n^{2-\varepsilon}$ that meets every sufficiently large square in the grid. Green records Ellenberg's grid question, brought to his attention by Zachary Chase, at the end of his discussion of Problem 20 in~\cite{green}.

A first lower bound comes from a single row, where isosceles triangles are exactly three-term arithmetic progressions. By Behrend's construction~\cite{behrend}, there is an absolute $C>0$ such that $F(n)\ge n\exp(-C\sqrt{\log n})$ for all sufficiently large $n$.

Charton, Ellenberg, Wagner and Williamson~\cite{patternboost} sampled points in the whole grid and deleted vertices of surviving triangles to obtain $F(n)\gg n/\sqrt{\log n}$. They also searched small grids computationally and suggested that a random greedy process might yield a linear lower bound.

Song~\cite{song} obtained the lower estimate below by sparse hypergraph colouring, and Croot, Mao, Pohoata, Sheffer and Yip~\cite{croot} proved the upper estimate. For an absolute $c_0>0$, their bounds give
\[
 n\sqrt{\frac{\log\log n}{\log n}} \ll F(n) \ll n^2\exp\!\left(-c_0\cdot\frac{\log n}{\log\log n}\right).
\]
Theorem~\ref{thm:main} removes the logarithmic loss in the lower bound by a deterministic construction.

The sieve of Croot et al.\ uses the fact that in an isosceles-free set, distinct pairs at a fixed distance cannot share a vertex. For some $c_1>0$, they bound the size of a grid set with all pairwise distances distinct by $O(n\exp(-c_1\log n/\log\log n))$. Thus for large $n$ an isosceles-free set as large as in Theorem~\ref{thm:main} must repeat a distance on disjoint pairs. For Riddell's question they announce that for some $\eta>0$ there are $m$-point planar sets whose isosceles-free subsets have size $O(m^{1-\eta})$. For the square grid, it remains open whether $F(n)/n$ is unbounded and whether there is a fixed $\delta>0$ such that $F(n)\ll n^{2-\delta}$.

Green also discusses the problem of avoiding isosceles right triangles of arbitrary orientation. Let $R(n)$ denote the largest size of a subset of $I_n^2$ containing no such triangle. In Gaussian integer coordinates, let $z$ be the right-angle vertex and choose the leg vector $w$ whose counterclockwise quarter-turn gives the other leg. The vertices are then $z,z+w,z+iw$. Pilatte~\cite{pilatte} showed that for some $\eta>0$, any subset of the grid avoiding these triangles has size $\ll n^2/(\log n)^{1+\eta}$. The stronger bound of Jaber, Liu, Lovett, Ostuni and Sawhney~\cite{corners} for axis-parallel corners applies here as well, since each such corner is an isosceles right triangle.

K\'arolyi and Solymosi~\cite{karolyi} used Gaussian integer digits to show that $R(n)\ge n^{\alpha-o(1)}$, where $\alpha=\log 281/\log(51\sqrt{2})=1.317865\ldots$. Our second theorem improves this exponent.

\begin{theorem}\label{thm:right}
For every $\varepsilon>0$ there is a constant $c_\varepsilon>0$ such that $R(n)\ge c_\varepsilon n^{2\log_5 3-\varepsilon}$ for all $n\ge1$.
\end{theorem}

Here $2\log_5 3=1.365212\ldots$. Modulo the Gaussian prime $2-i$, the equation for an isosceles right triangle becomes a midpoint equation. We choose digit vectors on a sphere to rule out nonconstant solutions.

\section{Proof of the linear lower bound}
\label{sec:proof}

For a prime $p$, write $\F_p$ for the field of residues modulo $p$, represented by the integers in $I_p$. We place scaled copies of smaller isosceles-free sets in distinct residue classes on one line modulo $p$. The first lemma forces any noncollinear isosceles triangle into one copy. The second lets us label the copies so that the same holds for equally spaced collinear triples.

\begin{lemma}\label{lem:line}
Let $m\ge1$ be an integer, let $p>4(m+1)^2$ be prime, and fix $t\in\F_p$. If the vertices of a noncollinear isosceles triangle in $I_{pm}^2$ lie on the line $y-(2m+1)x\equiv t\pmod p$, then they are congruent coordinatewise modulo $p$.
\end{lemma}

\begin{proof}
Let $a$ be the apex and $b,c$ the other vertices. Set $u=b-c$ and $v=b+c-2a$. Both are nonzero because $b\ne c$ and $a$ lies off the line through $b$ and $c$. The two equal sides give
\[
 u\cdot v=\lVert b-a\rVert_2^2-\lVert c-a\rVert_2^2=0.
\]
Every difference of two points on the prescribed line is a multiple of $(1,2m+1)$ modulo $p$. This holds for $u$ and also for $v=(b-a)+(c-a)$. Since $0<1+(2m+1)^2<p$, two nonzero multiples of this direction cannot be orthogonal modulo $p$. Thus at least one of $u,v$ is zero modulo $p$. Interchanging them if necessary, write $u=pw$, where $w=(w_1,w_2)$ is a nonzero integer vector. The coordinates of $u$ and $v$ have absolute value less than $2pm$, so $|w_1|,|w_2|<2m$.

If $w_1+(2m+1)w_2=0$, then $2m+1$ divides $w_1$, forcing $w_1=w_2=0$. Thus $w\cdot(1,2m+1)$ is a nonzero integer of absolute value less than $2m(2m+2)<p$. Dividing $u\cdot v=0$ by $p$ over the integers gives $w\cdot v=0$. Modulo $p$, the only multiple of $(1,2m+1)$ orthogonal to $w$ is zero, so $v\equiv0\pmod p$ as well. Now $b\equiv c$ and $b+c\equiv2a$ modulo $p$. Since $p$ is odd, $a\equiv b\equiv c\pmod p$.
\end{proof}

\begin{lemma}\label{lem:labels}
For every sufficiently large prime $p$ and every prime $q\ge p^{1/5}$, there is a labelling $f:\F_p\to\F_q$ such that any three-term arithmetic progression in $\F_p$ whose labels form an arithmetic progression in $\F_q$ is constant.
\end{lemma}

\begin{proof}
Put $B=\lceil p^{1/32}\rceil$. Write each $r\in I_p$ with $32$ base-$B$ digits, allowing leading zeros. Colour $r$ by the parities of its digits and the sum of their squares. This uses at most $2^{32}(1+32(B-1)^2)=O(p^{1/16})$ colours.

Let $r_1,r_2,r_3\in I_p$ have the same colour and satisfy $r_1+r_3\equiv2r_2\pmod p$. Corresponding digits of $r_1$ and $r_3$ have the same parity, so averaging them gives the digits of the integer $(r_1+r_3)/2$. This integer lies in $I_p$ and is congruent to $r_2$ modulo $p$, so it equals $r_2$. The digit vector of $r_2$ is therefore the midpoint of the other two. All three lie on the same sphere, since the colour fixes the sum of squared digits. The midpoint of two distinct points on a sphere lies strictly inside it, so $r_1=r_2=r_3$.

For labels, use the integers with $\lfloor\log_4 q\rfloor$ base-four digits in $\{0,1\}$, allowing leading zeros. All lie below $q/2$, so a midpoint equation modulo $q$ holds over the integers. In base four there are no carries, and comparing digits forces all three labels to agree. The number of labels is $2^{\lfloor\log_4 q\rfloor}\ge\sqrt q/2\ge p^{1/10}/2$, which exceeds the number of colours for large $p$. Assign a distinct label to each colour and let $f(r)$ be the label assigned to the colour of $r$. A progression whose labels form a progression then has a single colour and is therefore constant.
\end{proof}

\begin{proof}[Proof of Theorem~\ref{thm:main}]
For each $n\ge3$, let $p$ be the largest prime at most $n^{3/4}$ and set $m=\lfloor n/p\rfloor$. We shall choose an absolute constant $c>0$ and prove the following stronger claim by induction on $n$. For every $t\in\F_p$, the line $y-(2m+1)x\equiv t\pmod p$ contains an isosceles-free subset of $I_n^2$ of size at least $c(n+1)$. Prescribing the intercept $t$ will let us choose the label of each smaller copy.

Bertrand's postulate gives $p\asymp n^{3/4}$ and hence $m\asymp n^{1/4}$. For $n$ large enough that $3\le m<n$, let $q$ be the largest prime at most $m^{3/4}$. Again by Bertrand's postulate, $q\asymp m^{3/4}\asymp n^{3/16}$. Since $p^{1/5}\asymp n^{3/20}$ and $(m+1)^2\asymp n^{1/2}$, for all sufficiently large $n$ we have
\[
 p>4(m+1)^2,\qquad p>q\ge p^{1/5}.
\]
Choose $n_0\ge3$ so that these inequalities and $3\le m<n$ hold for every $n>n_0$ and $p$ is large enough for Lemma~\ref{lem:labels}. Set $c=1/(n_0+1)$. For $n\le n_0$ a single point suffices since $c(n+1)\le1$. When $n\ge3$, take $(0,t)$ with $t$ represented in $I_p\subseteq I_n$.

Fix $n>n_0$ and $t\in\F_p$, and assume the claim for smaller grid sizes. For each $r\in I_p$, let $v_r$ be the unique point of $I_p^2$ on the prescribed line with first coordinate $r$. Every point of $v_r+pI_m^2\subseteq I_{pm}^2$ reduces to $v_r$ modulo $p$. Thus these $p$ copies are disjoint and lie on the prescribed line in $I_n^2$.

Choose a labelling $f:\F_p\to\F_q$ from Lemma~\ref{lem:labels}. In the induction hypothesis at side length $m$, the prime is $q$, so the slope is $2\lfloor m/q\rfloor+1$. In copy $r$ we want every selected point $(x,y)$ to satisfy
\[
 y-\bigl(2\lfloor m/q\rfloor+1\bigr)x\equiv f(r)\pmod q.
 \numberthis\label{eq:intercept}
\]
Dilation by $p$ multiplies the intercept by $p$, and translation by $v_r$ adds a fixed residue. Since $p$ is invertible modulo $q$, exactly one starting intercept gives $f(r)$ in copy $r$. Apply the induction hypothesis in $I_m^2$ with that intercept. After dilation and translation, we obtain an isosceles-free set $A_r$ in copy $r$ with $|A_r|\ge c(m+1)$.

Let $A$ be the union of the sets $A_r$. By Lemma~\ref{lem:line}, the vertices of any noncollinear isosceles triangle in $A$ would be congruent modulo $p$ and hence lie in a single $A_r$.

\begin{samepage}
For a collinear isosceles triangle in $A$, let $r_1,r_2,r_3$ be the first coordinates modulo $p$, with $r_2$ corresponding to the apex. The apex is the midpoint, so $r_1+r_3=2r_2$ in $\F_p$. The left-hand side of \eqref{eq:intercept} is linear in $x,y$, so its values $f(r_1),f(r_2),f(r_3)$ also form a progression, now modulo $q$. Lemma~\ref{lem:labels} gives $r_1=r_2=r_3$, again putting all three vertices in one $A_r$. Since each $A_r$ is isosceles-free, so is $A$.
\par\end{samepage}

Finally, $m=\lfloor n/p\rfloor$ gives $n+1\le p(m+1)$. Summing over the disjoint sets gives the required bound
\[
 |A|\ge c\,p(m+1)\ge c(n+1). \qedhere
\]
\end{proof}

\section{Isosceles right triangles}
\label{sec:right}

\begin{proof}[Proof of Theorem~\ref{thm:right}]
Three distinct Gaussian integers form an isosceles right triangle if and only if they can be ordered as $a,b,c$ with $c-a=i(b-a)$. Put $\pi=2-i$. The residue field $\mathbb Z[i]/(\pi)$ is $\F_5$, and $i\equiv2\pmod{\pi}$. Thus the triangle equation reduces to the midpoint equation $a+c=2b$.

We choose our digits using Behrend's sphere construction~\cite{behrend,elsholtzpach}. For an integer $d\ge1$, partition $\{0,1,2\}^d$ according to the sum of the squares of the coordinates, and let $D$ be a largest class. There are at most $4d+1$ possible sums, so $|D|\ge3^d/(4d+1)$. If $u,v,w\in D$ satisfy $u+w=2v$ modulo $5$, the equality holds over the integers because each coordinate on either side lies between $0$ and $4$. The three vectors lie on the same sphere, and the midpoint of two distinct points on a sphere lies strictly inside it. Hence $u=v=w$.

For each integer $k\ge1$, form the set of vectors
\[
 W_k=\left\{\sum_{j=0}^{k-1}\pi^j u_j : u_0,\ldots,u_{k-1}\in D\right\}\subset\mathbb Z[i]^d.
\]
Reduction modulo $\pi$ recovers the first digit $u_0$. Subtracting it and dividing by $\pi$ recovers the remaining digits in turn, so $|W_k|=|D|^k$. Now suppose $a,b,c\in W_k$ satisfy $c-a=i(b-a)$ coordinatewise. Their first digits satisfy the midpoint equation in $D$, so all three first digits agree. Subtract the common digit from each vector and divide by $\pi$. Repeating this argument gives $a=b=c$.

To place $W_k$ in the plane, put $M=\lceil10|\pi|^k\rceil$ and encode a vector $w$ by the Gaussian integer
\[
 L(w)=w_1+Mw_2+\cdots+M^{d-1}w_d.
\]
If $L(e)=0$ for a Gaussian integer vector $e$ with $|e_j|<M$ for every $j$, then $M$ divides $e_1$, so $e_1=0$. Dividing by $M$ and repeating gives $e=0$.

Each coordinate of a vector in $W_k$ has modulus at most $2\sum_{j=0}^{k-1}|\pi|^j<2|\pi|^k$. Applying the preceding observation to differences shows that $L$ is injective on $W_k$. If three images satisfy the triangle equation, then $L(c-a-i(b-a))=0$. Each coordinate of $c-a-i(b-a)$ has modulus less than $8|\pi|^k<M$, so $c-a=i(b-a)$ in $W_k$. We have already shown that this forces $a=b=c$.

Since $M\le11|\pi|^k$ and $|\pi|=\sqrt5$, the $|D|^k$ images fit after translation in a grid of side at most $C_d5^{kd/2}$ for a constant $C_d$ depending only on $d$. For large $n$, choose the largest $k$ with $C_d5^{kd/2}\le n$. Then $5^{kd/2}\asymp_d n$, giving
\[
 R(n)\gg_d n^{(2/d)\log_5 |D|},\qquad
 \frac{2}{d}\log_5 |D|\ge 2\log_5 3-\frac{2\log_5(4d+1)}{d}.
\]
Choose $d$ large enough that $2\log_5(4d+1)/d\le\varepsilon$. This proves the stated bound for large $n$, and decreasing $c_\varepsilon$ covers the remaining grid sizes.
\end{proof}

\end{document}